\documentclass[11pt]{article}
\usepackage{geometry}                
\usepackage{graphicx}
\usepackage{amssymb}
\usepackage{amsmath}
\usepackage{epstopdf}
\usepackage{amsthm}
\usepackage{threeparttable}
\usepackage{booktabs}
\usepackage{setspace}
\usepackage[ruled,vlined]{algorithm2e}
\def\lj38{{\rm LJ}_{38}}

 \title{Metastability, Spectra, and Eigencurrents \\of the Lennard-Jones-38 Network}
 
 \author{Maria Cameron$^1$ }

\begin{document}
\maketitle

 \footnotetext[1]{University of Maryland, 
 Department of Mathematics, 
 College Park, MD  20742,
 {\tt cameron@math.umd.edu}}
 
\abstract{
We develop computational tools for spectral analysis of stochastic networks representing energy 
landscapes of atomic and molecular clusters.
Physical meaning  and some properties of eigenvalues, left and right eigenvectors, and eigencurrents are discussed.
We propose an approach to compute a collection of eigenpairs and corresponding eigencurrents describing the most important relaxation
processes taking place in the system on its way to the equilibrium.  
It is suitable for large and complex stochastic networks where  pairwise transition rates, 
given by the Arrhenius law, vary by orders of magnitude.
The proposed methodology is applied to the 
network representing the Lennard-Jones-38 cluster created by Wales's group. 
Its energy landscape has a double funnel structure with a deep and narrow face-centered cubic funnel and 
a shallower and wider icosahedral funnel. Contrary to the expectations, there is no appreciable
spectral gap separating the eigenvalue corresponding to the escape from the icosahedral funnel.
We  provide a detailed description of the 
escape process from the icosahedral funnel using the eigencurrent and
demonstrate a superexponential growth of the corresponding eigenvalue. 
The proposed spectral approach  is compared to the methodology of the Transition Path Theory.
Finally,  we discuss whether the Lennard-Jones-38 cluster is metastable from the points 
of view of a mathematician and a chemical physicist, and make a
connection with experimental works.
}
Modeling by means of stochastic networks has become widespread in chemical physics.
Efficient methods for the conversion of energy landscapes into stochastic networks have been developed.
Wales's group constructed a large number of networks representing energy landscapes 
of atomic and molecular clusters and proteins \cite{web,wales_book,wales_landscapes}. 
A number  of researchers 
developed efficient techniques for 
constructing continuous-time Markov chains representing coarse grained dynamics of biomolecules
and protein folding called Markov State Models \cite{schuette_thesis,swope,pande05,pande07,noe07,noe09,prinz,schuette11}.

Stochastic networks are attractive for analysis. On one hand, they are
more mathematically tractable than the original continuous systems. 
On the other hand, they are designed to  preserve important features
of the underlying systems.
Nevertheless,  analysis of large and complex networks still presents a challenge.
The numbers of states (vertices) and edges can be of the order of $10^n$, $n=3,4,5,6,...$,
and the pairwise transition rates may vary by tens of orders of magnitude. 
The study of such networks motivates development of new techniques and leads to a new level of understanding of 
complex systems.

In this work, we focus on developing a methodology for computing eigenvalues, eigenvectors, and eigencurrents 
for networks representing energy landscapes whose
pairwise transition rates are given by the Arrhenius law. The spectral decomposition of the generator matrix of a stochastic network
exposes all of the relaxation processes taking place in the system on its way to the equilibrium.  
The concept of metastability  and its relationship with the spectral structure has received a lot of attention 
in the scientific literature \cite{bovier2002,bovier1,bovier2,bovier3,gaveau,kurchan1,schuette03,schuette04}.
The assumption of the existence of a spectral gap separating a few smallest eigenvalues  allows us to 
simplify the dynamics of a complex stochastic system by 
combining its states into quasi-invariant subsets and considering transitions between them.
As the temperature tends to zero, 
the  eigenvectors can be approximated by the indicator functions 
of those quasi-invariant subsets, while the eigenvalues
give the exit rates from them \cite{bovier2002,bovier1}.

We have recently shown that, under certain nonrestrictive assumptions,
the zero-temperature asymptotics for  eigenvalues and eigenvectors of the generator matrices of such networks can be found recursively
starting from the low-lying group \cite{cam2}. 
This recursion can be realized in the form of an efficient algorithm which, in a nutshell, 
is a procedure for removing edges one-by-one 
 in a certain order from the minimum spanning tree 
for the given network \cite{cam2}.
Each asymptotic eigenpair is associated with a particular state of the stochastic network.
This state has the smallest potential value out of those where the asymptotic eigenvector is nonzero.

Here we upgrade the zero temperature asymptotic analysis with a continuation technique for computing a collection of 
eigenpairs of interest at a finite temperature range.  
It is based on the Rayleigh
Quotient iteration supplemented with some precautions. 
Each eigenpair defines a relaxation process from a perturbed probability distribution where the perturbation,
given by the left eigenvector, decays uniformly throughout the network, and the decay rate is given by the eigenvalue.
The quantitative description of the relaxation process for each eigenpair is given by the corresponding eigencurrent \cite{cve}.
We provide its physical interpretation and prove some of its properties. 
Eigencurrents can be  readily computed from the corresponding finite temperature eigenvectors.
At each temperature, the collection of states can be partitioned into those where the eigencurrent is emitted and those where it is absorbed.
The collection of edges connecting emitting and absorbing states 
can be viewed as the true transition state for the relaxation process. We will call it the emission-absorption cut.

It is instructive to compare the key concepts of the spectral analysis with those of the 
Transition Path Theory (TPT)  \cite{eve1,eve2,dtpt,luna}. We will discuss similarities and differences between the 
escape process quantified by the left and right eigenvectors, the eigenvalue, and the eigencurrent and the Transition Path Process \cite{luna, cve}
described by the committor, the transition rate, and the reactive current.

\begin{figure}
\centerline{
\includegraphics[width=0.5\textwidth]{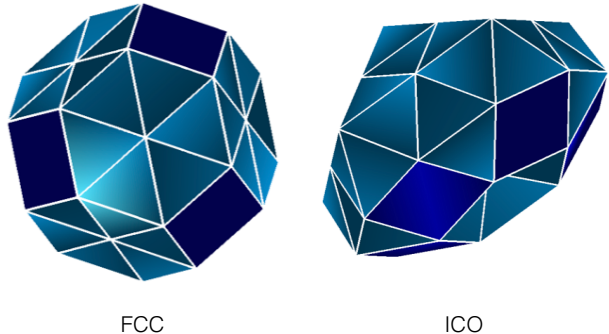}
}
\vspace{-5mm}
\caption{\label{fig:fico} 
Left: the global minimum of the $\lj38$ cluster, 
the face-centered cubic truncated octahedron with the point group $O_h$.
Right: the second lowest minimum of the $\lj38$ cluster, 
an incomplete icosahedron with the five-fold rotational symmetry (the point group $C_{5v}$). 
}
\end{figure}

We apply the proposed methodology to the Lennard-Jones-38  ($\lj38$) network with 71887 states (local minima) and 119853 edges
(transition states) created 
by Wales's group and available on the web \cite{wales_network}.
The $\lj38$ cluster has been profoundly studied by a number of researchers
and used as a test problem for various newly developed techniques 
\cite{wales38,wales_thermo,wales0,picciani,neirotti,voter,cam1,cve,cam2}.
The $\lj38$ cluster is remarkable due to the fact that it is the smallest cluster whose global minimum is achieved at a configuration
based on other than icosahedral packing \cite{wales-doye}.   
Its energy landscape   
has a double-funnel structure \cite{wales38,wales_landscapes,wales_book}.
The deep and narrow face-centered cubic (fcc) funnel ends at  the global
minimum, the face-centered 
cubic  truncated octahedron with the  point group $O_h$ (Fig. \ref{fig:fico}, left). We will denote it by FCC.
The wider and shallower icosahedral funnel is crowned by  the second lowest minimum, 
a five-fold rotationally symmetric incomplete icosahedron with the point group $C_{5v}$ (Fig. \ref{fig:fico}, right). 
We will denote it by ICO.
At  temperatures close to zero, the system spends most of the time at the fcc funnel, 
while at higher temperatures,
the icosahedral funnel with the greater configurational entropy becomes more 
preferable  \cite{wales_thermo,wales_book,frantsuzov}.
The $\lj38$ stochastic network is a feasible model for the low-temperature dynamics of the $\lj38$ cluster.

The double-funnel structure of the energy landscape might make us expect that the  
smallest nonzero eigenvalue of the generator matrix of the $\lj38$ network 
corresponds to the escape process from the icosahedral funnel.
However, this eigenvalue is buried under the number \cite{cam2} 245. Moreover, there is no appreciable 
 spectral gap separating it from the rest. 
 There are a few notable gaps separating eigenvalues corresponding to the escape processes from some
isolated  high-lying liquid-like states over extremely high potential barriers. 
The rest of the spectrum is pretty dense,
and  gaps appear in it  only at temperatures very close to zero. 
Roughly speaking, the spectrum is polluted by a large number of 
eigenvalues corresponding to the escapes from isolated high-lying  
liquid-like states, while these states are dynamically accessible from
icosahedral and fcc states, only if the waiting time is extremely long.
If we bound the observation time \cite{wales2014} which is equivalent to imposing a cap on the potential energy, 
the expected spectral gap separating the eigenvalue, responsible for the escape from the icosahedral funnel, appears. 

We have computed  a collection of finite temperature 
eigenpairs for the $\lj38$ network corresponding to the escape processes from the
largest quasi-invariant sets. 
The eigenvalue 
corresponding to the escape process from the icosahedral funnel exhibits a superexponential growth.
A quantitative description of this escape process at different values of temperature is provided by the corresponding eigencurrent.
The emission distribution of the eigencurrent together with the  distribution of
the eigencurrent in the emission-absorption cut explain the observed superexponential temperature  dependence of
the escape rate from the icosahedral funnel.

Finally, 
we discuss the issue of metastability of $\lj38$. The absence of the spectral gap renders this question nontrivial
from the mathematical point of view. On the other hand, the set of icosahedral configurations is certainly metastable 
from the point of view of a chemical physicist.
We undertake an attempt to make a connection between our results obtained for the $\lj38$ network and 
the experimental results where rare gas clusters such as argon, krypton and xenon have been studied using
electron diffraction \cite{farges,echt1,harris1,harris2,echt2,kovalenko1,kovalenko2} or x-rays \cite{kakar}.

The rest of the paper is organized as follows. Section \ref{sec:sig} 
is devoted to a physical interpretation  and some properties of eigenvalues, left and right eigenvectors, 
and eigencurrents of the generator matrices 
of networks with detailed balance. 
In Section \ref{sec:as}, we give an overview of zero-temperature asymptotic estimates for low-lying spectra of stochastic networks.
 Numerical algorithms are presented in Section \ref{sec:methods}. 
Section \ref{sec:lj38} is devoted to the application to the $\lj38$ network. 
The analyses of the $\lj38$ network  by means of the proposed spectral approach and the Transition Path Theory  \cite{eve1,eve2,dtpt,luna,cve}
are compared in Section \ref{sec:tpt}.  Metastability issues and a connection with experimental works are 
discussed in Section \ref{sec:met}.
Conclusions are summarized  in Section \ref{sec:con}.

\section{\label{sec:sig}Significance of spectral analysis}
In this Section, we discuss what can we learn from the spectral decomposition 
of the generator matrix describing the dynamics of a stochastic network.

\subsection{\label{sec:eig}The eigenstructure of networks with detailed balance}
In this Section, we provide some background and introduce some notations.
Let $G(S,E,L)$ be a network, where  $S$  and $E$ are its sets of states and edges respectively.
We assume that the number of states is finite and denote it by $N$.
The dynamics of this network  is described by the generator matrix $L =\{L_{ij}\}_{i,j=1}^N$. If states $i$ and $j$ ($i\neq j$)
are connected by an edge, then $L_{ij}$ 
is the pairwise transition rate from $i$ to  $j$, otherwise $L_{ij}=0$. 
 The diagonal entries $L_{ii}$ are defined so that the row sums of the matrix $L$ are zeros. 
 This fact implies that $L$ has a zero eigenvalue
$\lambda_0=0$ whose right eigenvector can be chosen to be $e:=[1,1,\ldots,1]^T$, i.e., $Le=0$. 
The corresponding left eigenvector can be chosen to be  the equilibrium probability distribution $\pi=[\pi_1,\pi_2,\ldots,\pi_N]$: 
$$\pi^TL=0,\quad\sum_{i=1}^N\pi_i= 1.$$
The generator matrices of networks representing energy landscapes possess the detailed balance property,
 i.e., $\pi_iL_{ij}=\pi_jL_{ji}$, which means that on average, there is the same numbers of transitions from 
 $i$ to $j$ and from $j$  to $i$ per unit time. 
 The detailed balance implies that the matrix $L$ can be decomposed as 
 $$L=P^{-1}Q$$
 where $P$ is the diagonal matrix 
 $$P={\rm diag}\{\pi_1,\pi_2,\ldots,\pi_N\},$$
 and $Q$ is symmetric. Hence $L$ is similar to a symmetric matrix:
 \begin{equation}
 \label{lsym}
 L_{sym}:=P^{1/2}LP^{-1/2}=P^{-1/2}QP^{-1/2}.
 \end{equation}
Using Eq. \eqref{lsym} and the fact that $L_{ii}=-\sum_{j\neq i}L_{ij}$ where $L_{ij}\ge0$, one can show 
that eigenvalues of $L$ are real and  nonpositive. 
Furthermore, we assume that the network is connected.  Hence $L$ is irreducible and 
the zero eigenvalue has algebraic multiplicity one.
We will denote the nonzero eigenvalues of $L$ by $-\lambda_k$, $k=1,2,\ldots,N-1$, 
and order them so that
$$0<\lambda_1\le\lambda_2\le\ldots\le\lambda_{N-1}.$$

The eigendecomposition of a generator matrix $L$
leads to a nice representation of the time evolution of the probability
distribution. The matrix $L$ can be written as 
\begin{equation}
\label{ed}
L=\Phi\Lambda\Phi^TP,
\end{equation}
where $\Phi=[\phi^0,\phi^1,\ldots,\phi^{N-1}]$ 
is a matrix of right eigenvectors, normalized so that $P^{1/2}\Phi$ is orthogonal.
Note that $P^{1/2}\Phi$ is the matrix of eigenvectors of $L_{sym}$.
Therefore, 
\begin{align}
\|P^{1/2}\phi^k\|^2&=\sum_{j=1}^N\pi_j|\phi^k_j|^2=1,\label{o1}\\
\sum_{j=1}^N\pi_j\phi^k_j\phi^l_j& = 0,\quad k\neq l,\label{o2}
\end{align}
and in particular,
$\phi^0=e$ and  $P\phi^0=\pi$.

The probability distribution evolves according to the 
 Forward Kolmogorov (a. k. a. the Fokker-Planck) equation
$$\frac{dp}{dt}=L^Tp,\quad p(0)=p_0,$$ 
where $p_0$ is the initial distribution.
Using Eqs. \eqref{ed}--\eqref{o2} we get
\begin{equation}
\label{eq1}
p(t) = \sum_{k=0}^{N-1}c_ke^{-\lambda_kt}P\phi^k, ~~{\rm where}~~c_k=(\phi^k)^Tp_0.
\end{equation}
Note that $c_0=e^Tp_0=\sum_{j=0}^N(p_0)_j=1$.
Eq. \eqref{eq1} 
shows that the $k$-th  eigen-component of $p(t)$ remains significant 
only on the time interval $O(\lambda_k^{-1})$.
Eventually, all eigen-components, except for the zeroth, decay.  
Hence $p(t)\rightarrow\pi$ as $t\rightarrow\infty$. 
\subsection{\label{sec:eve}Interpretation of eigenvalues and eigenvectors}
For $k>0$, the left and right  $k$-th eigenvectors of $L$, $P\phi^{k}$ and $\phi^{k}$ 
respectively, can be understood from a recipe for preparing the initial probability distribution so that  only the coefficients
$c_0$ and $c_k$ in Eq. \eqref{eq1} are nonzero.
Imagine $n\gg1$ particles distributed  in the stochastic network according to the equilibrium distribution $\pi$.
Since $(\phi^0)^TP\phi^k=\sum_{j=1}^N\pi_j\phi_j^k=0$ for any $k>0$, 
the set of states $S$ can be divided into two parts: 
\begin{align}
S_{+}^k &:= \{i\in S~:~\phi^k_i\ge0\}\quad{\rm and}\notag\\
S^k_{-}&:= \{i\in S~:~\phi^k_i<0\}.\label{pm}
\end{align}
In order create a component of the initial distribution
parallel to the left eigenvector $P\phi^k$, we pick some $\alpha$ satisfying   $\alpha\phi^k_j\le1$ for all $j\in S^k_{-}$, 
remove $\alpha n \pi_j|\phi^k_j|$ particles from each state $j\in S^k_{-}$ and 
distribute them in  $S^k_{+}$ so that each state $j\in S^k_{+}$ gets $\alpha n\pi_j\phi^k_j$ particles. 
The resulting distribution is 
$$p_0 = \pi+\alpha P\phi^k.$$ 
Therefore, the left eigenvector $\alpha P\phi^k$ is a perturbation of the  
equilibrium distribution decaying uniformly across the network with the rate 
given by the corresponding eigenvalue 
$\lambda_k$. (The keyword here is ``uniformly".) The corresponding right eigenvector $\alpha\phi^k$ 
for each state $j\in S$ gives the factor by which it is over- or underpopulated with respect to the equilibrium distribution $\pi$, as
$$\alpha\phi^k_j \equiv \frac{ (p_0)_j - \pi_j}{\pi_j}.$$
 Therefore, $\phi^k$ is,
 in essence, a fuzzy signed indicator function of the perturbation $P\phi^k$.

\subsection{Eigencurrents}
Probability currents are important tools for the quantitative description 
of transition processes in the system \cite{kurchan6}.
For example, the reactive current is one of the key concepts of the Transition Path Theory \cite{eve1,eve2,dtpt,cve}.
In the context of spectral analysis, E. Vanden-Eijnden proposed to consider  eigencurrents \cite{eve_chapter,cve}.
While eigenvectors determine the perturbations 
to the equilibrium distribution decaying with the rates given by the corresponding eigenvalues, eigencurrents
give a quantitative description of the escape process from these perturbed distributions.
 
Eigencurrents are defined as follows. The time derivative of the $i$-th component of the probability distribution $p(t)$
is
\begin{equation}
\label{eq2}
\frac{dp_i}{dt}=\sum_{j=1}^NL_{ji}p_j=\sum_{j\neq i}(L_{ji}p_j-L_{ij}p_i).
\end{equation}
Plugging in  expressions for $p_i$ and $p_j$ from Eq. \eqref{eq1} into Eq. \eqref{eq2}  
and using the detailed balance property $\pi_iL_{ij}=\pi_jL_{ji}$
we obtain
\begin{equation}
\label{eq3}
\frac{dp_i}{dt}=\sum_{k=0}^{N-1}c_ke^{-\lambda_kt}\sum_{j\neq i}\pi_iL_{ij}[\phi^k_j-\phi^k_i].
\end{equation}
The collection of numbers
\begin{equation}
\label{eq4}
F^k_{ij}:=\pi_iL_{ij}e^{-\lambda_kt}[\phi^k_i-\phi^k_j]
\end{equation}
is called the eigencurrent associated with the $k$-th eigenpair. 
(Here we have switched the sign and incorporated the factor $e^{-\lambda_kt}$ into the definition 
in comparison with the one  in Ref. \cite{cve} in order to make its physical sense more transparent.) 
In terms of the eigencurrent $F^k_{ij}$ Eq. \eqref{eq3} can be rewritten as 
\begin{equation}
\label{eq5}
\frac{dp_i}{dt}=-\sum_{k=0}^{N-1}c_k\sum_{j\neq i}F^k_{ij}.
\end{equation}
Hence, the eigencurrent $F^k_{ij}$ is  
the net probability current of the $k$-th perturbation $P\phi^k$ along the edge $(i, j)$ per unit time.
In other words, if the system is originally distributed according to $p_0=\pi+\alpha P\phi_k$ 
then the current $\alpha F^k_{ij}$ gives the difference of the
average numbers of transitions from $i$ to $j$ and from $j$ to $i$ per unit time.

The eigencurrent $F^k_{ij}$ can be compared to the reactive current $F^{R}_{ij}$ in Ref. \cite{cve} 
(which is the same as the quantity $f_{ij}^{AB}-f_{ji}^{AB}$ in Ref. \cite{dtpt}). The reactive current $F^R_{ij}$ 
 is conserved at every node of the network  \cite{dtpt} except for the specially 
designated subsets of source states $A$ and sink states $B$.
Contrary to it, the eigencurrent is either emitted or absorbed 
in all states $i$ where $\phi^k_i\neq 0$. Indeed, for any $i\in S$ we have 
\begin{align}
\sum_{j\neq i}F^k_{ij}=&\pi_ie^{-\lambda_kt}\left[\phi_i^k\sum_{j\neq i}L_{ij}-\sum_{j\neq i}L_{ij}\phi_j^k\right]  \notag\\
=&\pi_ie^{-\lambda_kt}\left[-\sum_{j=1}^NL_{ij}\phi_j^k\right]=e^{-\lambda_kt}\lambda_k\pi_i\phi_i^k. \label{eq6}
\end{align}
Therefore, every state $i$ emits or absorbs (depending on the sign of $\phi^k_i$) 
$e^{-\lambda_kt}\lambda_k\pi_i|\phi_i^k|$ units of the eigencurrent $F^k_{ij}$ per unit time.

Let us partition the set of states $S$ 
into  $S^k_{+}$ and $S^k_{-}$ (see Eq. \eqref{pm}). 
The corresponding cut-set (a. k. a. cut) consists of all edges $(i, j)$ where $\phi_i^k\ge0$ and $\phi_j^k<0$.
We will call this cut the emission-absorption cut as
it separates the states where the eigencurrent $F^k_{ij}$ is not absorbed from those where it is absorbed.
It  can be compared to the isocommittor cut \cite{cve} corresponding to the committor value $q=0.5$.  

Now we will show that 
for every fixed time $t$, the eigencurrent over the emission-absorption cut is maximal among all possible cuts of the network,
and it is equal to the total  eigencurrent emitted by the states $S^k_{+}$ per unit time at time $t$, i.e.,
\begin{align}
&\max_{S',S'':~S=S'\dot{\cup}S''}\sum_{i\in S',j\in S''}F^k_{ij}\notag\\
=&\sum_{i\in S^k_{+},j\in S^k_{-}}F^k_{ij}=e^{-\lambda_kt}\lambda_k\sum_{i\in S^k_{+}}\pi\phi^k_i\label{eq7}
\end{align}
(The symbol $\dot{\cup}$ denotes the disjoint union.)
Since Eq. \eqref{eq4} implies that $F^k_{ij}=-F_{ji}^k$,
for any subset $S'\subset S$ we have
\begin{equation*}
\sum_{i\in S',j\in S'}F^k_{ij}=0.
\end{equation*}
Therefore, the eigencurrent over the cut $(S',S'')$ is
\begin{align*}
 \sum_{i\in S',~j\in S''}F_{ij}^k &= \sum_{i\in S',~j\in S}F_{ij}^k  - \sum_{i\in S',~j\in S'}F_{ij}^k \notag \\
 =&e^{-\lambda_kt}\lambda_k\sum_{i\in S'}\pi_i\phi^k_i,
 \end{align*} 
 i.e., it is the total eigencurrent emitted in $S'$ per unit time. 
The maximum of the last sum is achieved if $S'=S^k_{+}$, i.e., if $S'$ consists of all non-absorbing states. 

The cut-set  of the emission-absorption cut is the true transition state 
of the relaxation process from the perturbed distribution $\pi+\alpha P\phi^k$.


\section{\label{sec:as}Asymptotic estimates for eigenvalues and eigenvectors}
Up to now, we have not made any assumptions about the 
particular form of the pairwise transition rates $L_{ij}$.
As a result, we have not offer any specific recipe for finding the set of eigenpairs (or its subset of interest).
This can be readily done  for small networks using direct methods. 
However, for large networks 
with  complex underlying potential energy landscapes and 
pairwise transition rates varying
by orders of magnitude, 
computing low-lying eigenpairs  is a challenging problem.  
Therefore, a priori estimates for eigenvalues and eigenvectors 
for such networks are of great importance.

In 1970's, A. Wentzell considered the spectral problem for stochastic 
networks without  the detailed balance but with pairwise transition 
rates of the form $e^{-U_{ij}/2\epsilon^2}$, 
where $\epsilon$ was a small parameter \cite{wentzell1,wentzell2} (see also Ref. \cite{f-w}, Chapter 6).
Such networks come, e. g.,  from continuous systems evolving according 
to non-gradient stochastic differential equations with small white noise \cite{f-w}.
Wentzell obtained asymptotic estimates  for eigenvalues up to the logarithmic order of the form $\lambda_k\asymp e^{-\Delta_k/\epsilon}$ where 
the numbers $\Delta_k$ are the  solutions of certain optimization problems over the so called $W$-graphs \cite{wentzell1,wentzell2,f-w}.

In 2000's, A. Bovier and collaborators developed sharp estimates for low lying  eigenvalues 
and the corresponding eigenvectors   for stochastic networks  with detailed balance
using the conceptual apparatus of the classic potential theory  \cite{bovier2002,bovier1}.
They also extended their results for continuous systems \cite{bovier2,bovier3}.
Their sharp estimates are obtained under the assumption of the existence of the hierarchy of sets of so called metastable points.
The estimates for eigenvalues are expressed in terms of expected exit times from certain subsets of states (see below) 
while  right eigenvectors are approximated by capacitors (a. k. a.  committors) for certain pairs of source and sink subsets of states.

The sharp estimates of Bovier et al are of great theoretical value. 
However, in order to implement them in practice, one actually needs to 
construct the hierarchy of  sets of metastable points. 
A recipe for doing this for networks representing energy landscapes was
proposed in Ref. \cite{cam2}. 
The nonzero off-diagonal entries of the generator matrix $L$ for such networks are 
of the form \cite{wales0}
$$L_{ij} = \frac{a_{ij}}{a_i}e^{-(V_{ij}-V_i)/T},$$
where $a_{i}$ and $a_{ij}$ are the prefactors 
for the local minimum $i$ and the saddle $(i,j)$ separating the local minima $i$ and $j$ respectively, 
$V_{i}$ and $V_{ij}$ are the values of the potential at  $i$ and  $(i,j)$ respectively, and $T$ is the temperature. 
This recipe involves the construction of the set of the optimal $W$-graphs and 
also relies on Bovier's et al estimates for  eigenvectors\cite{bovier2002,bovier1} and Wentzell's estimates for  eigenvalues\cite{wentzell1,wentzell2,f-w}.

In the graph theory, a tree is a connected graph with no cycles; 
a forest is a collection of trees;
a spanning tree for a given undirected graph is a tree with the same set of vertices; 
and a minimum spanning tree for a weighted undirected graph is a spanning tree 
such that the sum of the weights of its edges is minimal possible \cite{amo}. 
A $W$-graph with $k$ sinks is a directed forest
such that in each of its $k$ connected components (subtrees)
all directed paths lead to a unique vertex called sink.
The optimal $W$-graph with one sink is the minimum spanning tree where the cost of each edge $(i,j)$ is $V_{ij}$,
and the directions of the 
edges are chosen so that all directed paths lead to the state with the smallest value of the potential.
We will call it sink $s^{\ast}_0$.

We have proven  \cite{cam2} that for the networks representing energy landscapes, 
the optimal $W$-graphs  have a nested property
which allows us to construct them recursively from top to bottom by solving a 
simple optimization problem at each step.
At step $k$, we add one sink $s^{\ast}_k$ and remove one edge $(p^{\ast}_k,q^{\ast}_k)$.
Then the asymptotics for the $k$th eigenvalue is $\lambda_k\asymp e^{-\Delta_k/T}$, 
where $\Delta_k:=V_{p^{\ast}_kq^{\ast}_k}-V_{s^{\ast}_k}$.
The asymptotic eigenvector $\phi_k$ is the indicator function of the set of states
of the connected component of the optimal $W$-graph at step $k$ containing the sink $s^{\ast}_k$.
The sinks fit the definition of metastable points\cite{bovier2002,bovier1} if the temperature is low enough. 

\section{\label{sec:methods}Computing eigenvalues and eigenvectors}
In this Section, we propose a two-step approach for computing eigenvalues and eigenvectors.
Step one is the algorithm introduced in Ref. \cite{cam2} for finding zero-temperature asymptotics. 
Here we will give its brief description.
Step two is computing eigenpairs 
of interest at finite temperatures using the  asymptotic eigenvectors as initial guesses.

\subsection{\label{sec:case}Computing asymptotic eigenpairs}
For simplicity  we make the following genericness assumption: 
all numbers $V_i$, $V_{ij}$ and the differences $V_{pq}-V_i$ are different.
The cases where this assumption is violated will be addressed in future work. 
Note that the genericness assumption implies that 
 the minimum spanning tree is unique. 

The preprocessing for the algorithm \cite{cam2} is computation of the minimum spanning tree $\mathcal{T}^{\ast}$ 
for the given network, where 
the cost of  the edge $(i,j)$ is $V_{ij}$, the value of the potential at the saddle $(i,j)$.
In other words, one needs to 
extract a collection of $N-1$ edges
such that all $N$ states are connected by these edges, and the sum of the weights of these edges is minimal possible.
The rest of the edges are removed from the network. 
There are several greedy algorithms for performing this task \cite{amo}. 
We have used Kruskal's algorithm \cite{kruskal}.
Note that all edges constituting the minimum spanning tree are such that for any pair of states $i$ and $j$,
the maximal cost $V_{pq}$ along the 
unique path $w^{\ast}(i,j)$ in $\mathcal{T}^{\ast}$ connecting $i$ and $j$
is minimal possible \cite{amo}.

Central to the algorithm \cite{cam2} 
are the barrier function $u(i)$ and the escape function $v(i)$, $i\in S$.
For a given set of sinks  $s^{\ast}_j$, $j=0,1,\ldots,k$,  
\begin{align}
u(i) &= \min_{0\le j\le k}\max_{(p,q)\in w^{\ast}(i,s^{\ast}_j)}V_{pq},\label{u}\\
v(i) &= u(i)-V_i,\label{v}
\end{align}
where $w^{\ast}(i,s^{\ast}_j)$ is the unique path in the minimum spanning tree $\mathcal{T}^{\ast}$ connecting
states $i$ and $s^{\ast}_j$.
Hence, $u(i)$ is the height of the  lowest possible highest saddle  separating state $i$ from the set of  sinks 
$\{s^{\ast}_j\}_{j=0}^k$, while 
$v(i)$ is the minimal possible maximal potential barrier to overcome in order to escape from $i$ to one of the sinks.
\begin{figure}
\centerline{
\includegraphics[width=0.6\textwidth]{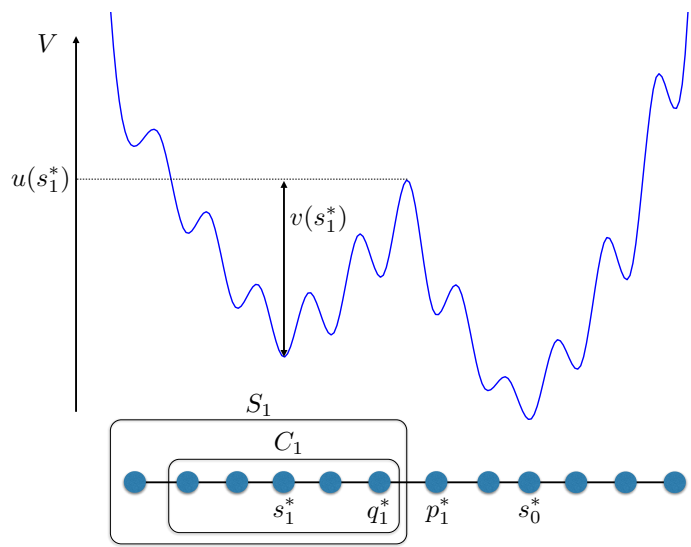}
}
\caption{\label{fig:pic} 
A chain-of-states example. In this case, the minimum spanning tree is the graph itself.
The sink $s^{\ast}_0$ is the state with the minimal value of the potential.  
The sink $s_1^{\ast}$ is the state for which the escape function $v$ computed with respect to 
$s^{\ast}_0$ is maximal. The cutting edge $(p^{\ast}_1,q^{\ast}_1)$ is the edge
with the maximal potential value in the unique path connecting $s^{\ast}_0$ and $s^{\ast}_1$.
Remove the edge $(p^{\ast}_1,q^{\ast}_1)$.
The quasi-invariant set $S_1$ is the set of states in the connected component of the resulting forest 
containing the new sink $s^{\ast}_1$. The Freidlin's cycle $C_1$ is the subset of $S_1$ consisting of states 
separated from $s_1^{\ast}$ by lower potential barriers than $V_{p^{\ast}_1q^{\ast}_1}$.
}
\end{figure}

\begin{figure}
\centerline{
\includegraphics[width=0.6\textwidth]{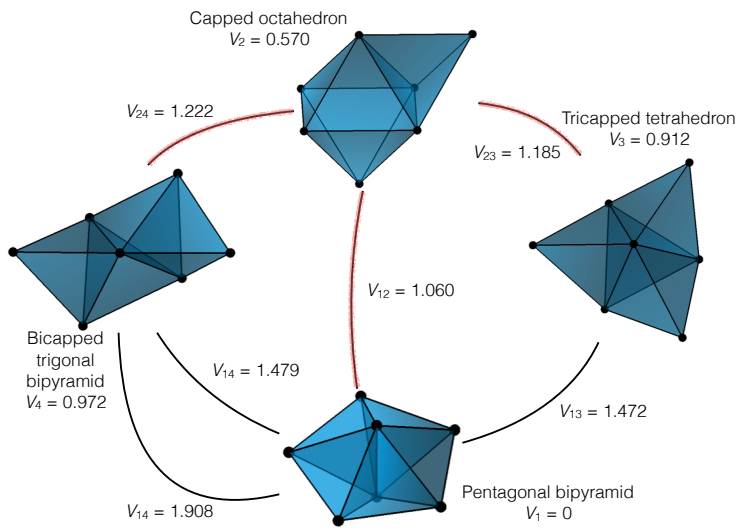}
}
\caption{\label{fig:lj7_network} 
A network associated with the energy landscape of the LJ$_7$ \cite{LJ7_tsai,LJ7_miller}. 
The values of $V_i$ and $V_{ij}$, $i,j\in\{1,2,3,4\}$ are given with respect to
the global potential minimum (the pentagonal bipyramid, $V = -16.505$).
The edges constituting the minimum spanning tree are marked red.
}
\end{figure}

\begin{figure}
\centerline{
\includegraphics[width=0.6\textwidth]{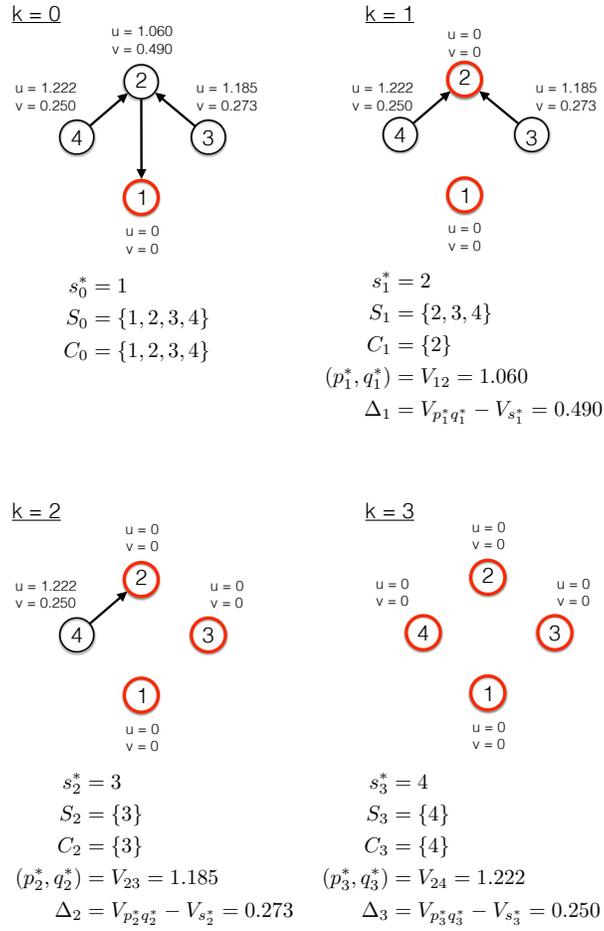}
}
\caption{\label{fig:lj7_algorithm} 
Initialization ($k=0$) and three iterations ($k=1,2,3$) of the algorithm for computing zero-temperature asymptotics for eigenvalues and eigenvectors
applied to the stochastic network in Fig. \ref{fig:lj7_network} associated with LJ$_7$.
}
\end{figure}

The algorithm \cite{cam2} proceeds as follows.
The initial sink is $s^{\ast}_0=\arg\min_{i\in S} V_i$.
The functions $u$ and $v$ are computed with respect to this sink. 
 The initial forest $\mathcal{T}_0^{\ast}$ is the minimum spanning tree 
$\mathcal{T}^{\ast}$. Set  $S_0=C_0=S$.
Then for $k=1,2,\ldots$
\begin{enumerate}
\item
Pick the state with the largest value of the escape function $v$ and make it the next sink $s^{\ast}_{k}$.
\item 
Find the edge with the maximal value $V_{pq}$ in the path connecting the new sink $s^{\ast}_{k}$ with one 
of the existing sinks $s^{\ast}_j$, $0\le j\le k-1$, denote it by $(p^{\ast}_k,q^{\ast}_k)$ and remove it,
i.e., set $\mathcal{T}^{\ast}_{k}=\mathcal{T}^{\ast}_{k-1}\backslash(p^{\ast}_k,q^{\ast}_k)$.
\item
Set $\Delta_k:=V_{p^{\ast}_kq^{\ast}_k}-V_{s^{\ast}_{k}}$. Set $S_{k}$ to be the set of states in the connected component 
of the forest $\mathcal{T}^{\ast}_{k}$ containing the new sink $s^{\ast}_{k}$.
\item
Update the values of $u$ and $v$ in  the connected component 
of the forest $\mathcal{T}^{\ast}_{k}$ containing the new sink $s^{\ast}_{k}$.
Denote by $C_k$ the set of states where the values of $u$ and $v$ 
have actually changed. 
\end{enumerate}
The first step of this algorithm for a chain-of-states example is illustrated in Fig. \ref{fig:pic}.
An application of this algorithm to the network in Fig. \ref{fig:lj7_network} 
associated with the Lennard-Jones cluster of 7 atoms (LJ$_7$) \cite{LJ7_tsai,LJ7_miller}  is worked out in Fig. \ref{fig:lj7_algorithm}. 

To summarize, for each $k=1,2,\ldots,N-1$ we get
\begin{itemize}
\item the quasi-invariant set $S_k$;
\item the sink $s^{\ast}_k$, which is the deepest minimum in the quasi-invariant set $S_k$;
\item the cutting edge $(p^{\ast}_k,q^{\ast}_k)$, which corresponds to the lowest possible highest saddle separating the sink $s^{\ast}_k$
from the previously chosen sinks $s^{\ast}_j$, $j=0,1,\ldots,k-1$;
\item the number $\Delta_k=V_{p^{\ast}_kq^{\ast}_k}-V_{s^{\ast}_{k}}$;
\item Freidlin's cycle\cite{freidlin-cycles,freidlin-physicad,f-w, cam1} $C_k\subseteq S_k$  which is the set of states separated from the sink $s^{\ast}_k$ 
by smaller potential barriers than $V_{p^{\ast}_kq^{\ast}_k}$.
It is the largest Freidlin's cycle \cite{cam2}  
containing the new sink and not containing any state with a smaller potential value.
Its significance will be explained in Section \ref{sec:met}.
\end{itemize}
The asymptotics for the eigenvalues and  the right eigenvectors are given by
\begin{align}
\lambda_k&\asymp e^{-\Delta_k/T},\label{lam}\\
\phi^k_i&\approx \begin{cases}1,&i\in S_k\\0,& i\notin S_k.\end{cases}\label{vec}
\end{align}
Note the ordering: $\Delta_1<\Delta_2<\ldots$. Hence the eigenpairs are computed in the 
increasing order of the corresponding eigenvalues. Therefore, one can stop the for-cycle whenever all 
eigenpairs of interest are computed.

Finally, we remark on what happens if the genericness assumption is violated.
If $\Delta_k=\Delta_{k+1}$ for some $k$ than two cases are possible. 
If the corresponding cutting edges belong to different connected components of the forest $\mathcal{T}^{\ast}_k$
then no error occurs.
If the cutting edges belong to the same connected component, then the asymptotic eigenvectors $\phi^k$ and $\phi^{k+1}$
might be determined wrongly, but the error does not propagate for larger $k$ in the for-cycle.
In other words, the possible effect from the violation of the genericness assumption is local.

\subsection{\label{sec:rq}Finite temperature continuation}
Difficulties in computing eigenpairs at finite (but still low) temperatures  are associated with the facts that $(i)$
all eigenvalues tend to zero as $T\rightarrow 0$, 
and the ones of interest can be extremely small at the desired range of temperatures;
$(ii)$ the eigenvalues can cross as the temperature increases; $(ii)$ the right eigenvectors can be nearly parallel.
Experimenting with a few continuation techniques we have found out
that the Rayleigh Quotient iteration (see e.g. Refs. \cite{demmel,trefethen}) 
with the asymptotic estimates  \eqref{vec} for the eigenvectors as initial guesses
is a feasible and  robust approach as soon as some precautions are taken.

The initial approximations given by Eq. \eqref{vec} are sufficiently accurate so that
 the Rayleigh Quotient iteration converges rapidly (cubically).
However, occasionally,  convergence to a wrong eigenpair  may occur. 
For example, this can happen for eigenpairs  number $k$ and $l$ 
such that $(i)$ the sets $S_k$ and $S_l$ are both large, $(ii)$ $S_k\subset S_l$, $(iii)$ the number of states in $S_l\backslash S_k$ is 
small in comparison with the number of states in $S_k$, and $(iv)$ the potential values at sinks $s_k^{\ast}$ and $s^{\ast}_l$ are close.
Fortunately, the convergence to a wrong eigenpair is easy to detect using the orthogonality relationship  \eqref{o2}. 

Now we provide some technical details of computing finite temperature eigenpairs.
Instead of $L$, we work with the matrix $L_{sym}$ given by Eq. \eqref{lsym}.
Its eigenvectors $\psi^k$ are orthonormal and relate 
to the right eigenvectors $\phi^k$ of $L$ via $\psi^k=P^{1/2}\phi^k$.
To compute the $k$th eigenpair of $L_{sym}$ at temperature $T$, 
we set the initial guess for the eigenvector $\psi^k$ as
\begin{equation*}
\left(\psi_i^k\right)_0(T)=\begin{cases}\pi^{1/2}_i(T),&i\in S_k,\\ 0,&i\notin S_k,\end{cases}
\end{equation*}
run the Rayleigh Quotient iteration and obtain an eigenpair $(\lambda,\xi)$. 
Then we check whether $(\lambda,\xi)$ is a correct eigenpair by calculating the dot product $\xi^T\left(\psi^k\right)_0$.
Since eigenvectors of $L_{sym}$ corresponding to distinct eigenvalues are orthogonal, this inner product is
close to zero if $(\lambda,\xi)$ is a wrong eigenpair. Otherwise, it is close to one.
If $(\lambda,\xi)$ is the desired eigenpair, we set $\lambda_k(T)=\lambda$ and $\phi^k(T)=P^{-1/2}(T)\xi$. 
If $(\lambda,\xi)$ is a wrong eigenpair, we repeat the Rayleigh Quotient 
iteration with a more accurate initial guess for the eigenvector (and/or eigenvalue) obtained from
already computed eigenpairs for the same $k$ but different values of temperature.


\section{\label{sec:lj38}Application to the LJ$_\mathbf{38}$ network}
The $\lj38$ dataset created by Wales's group \cite{wales_network} contains 100000 local minima and 138888 transition states.
Its largest connected component consists of 71887 local minima (states) and 
119853 transition states (edges) and contains the two deepest minima FCC and ICO.
We will analyze this connected component and refer to it as the $\lj38$ network.
Local minima, other than FCC and ICO, will be referred to by their indices in the list \cite{wales_network}. 
Here are the  indices of some important minima: FCC has index 1, ICO has index 7, the third lowest minimum has index 16, 
the lowest possible highest saddle separating ICO and FCC
is the transition state between minima 342 and 354.

\subsection{\label{sec:set}Settings and thermodynamics} 
Following Ref. \cite{wales0} we define pairwise transition rates $L_{ij}$ using the harmonic approximation.
If there is a transition state in the database \cite{wales_network} separating local minima $i$ and $j$ then
\begin{equation}
\label{har}
L_{ij}(T)=\frac{O_i\bar{\nu}_i^{\kappa}}{O_{ij}\bar{\nu}_i^{\kappa-1}}e^{-(V_{ij}-V_i)/T},
\end{equation}
where $O_i$, $V_i$, and $\bar{\nu}_i$ are, respectively, the point group order, the value of the potential energy,
and the geometric mean vibrational frequency for minimum $i$;
$O_{ij}$, $V_{ij}$, and $\bar{\nu}_{ij}$
are the same parameters for the saddle $(i,j)$; 
and $\kappa = 3\cdot38-6 = 108$ is
the number of vibrational degrees of freedom. Otherwise, if there is no such a transition state (i.e., there is no edge $(i,j)$), $L_{ij}=0$.
The equilibrium probability distribution is given by
\begin{equation}
\label{pi}
\pi_i(T)=\frac{1}{Z(T)}\frac{e^{-V_i/T}}{O_i\bar{\nu}_i^{\kappa}},\quad 
Z(T)=\sum_{i\in S}\frac{e^{-V_i/T}}{O_i\bar{\nu}_i^{\kappa}}.
\end{equation}

To choose the range of temperatures at which we will investigate the $\lj38$ network,
we use its critical temperature
of the solid-solid phase transition as a reference.
It has been established \cite{frantsuzov,wales_thermo,wales_book} 
that the solid-solid phase transition, when  fcc structures give place to  icosahedral packings, 
occurs in the $\lj38$ cluster at $T=0.12$. The outer layer starts to melt, i.e., icosahedral states become less populated than liquid-like states,
at $T=0.18$. These critical temperatures were obtained in the continuous setting
using either Monte-Carlo simulations \cite{frantsuzov} or the superposition method for 
determination of the density of states and an anharmonic approximation for the pairwise transition rates \cite{wales_thermo,wales_anhar}.
However, if the continuous $\lj38$ cluster is modeled by the described above network 
with the pairwise rates  given by Eq. \eqref{har}, the critical temperatures are shifted upward. 
One can determine them, e.g.,  by plotting the formally computed \cite{wales_thermo,wales_book} 
heat capacity at constant volume: the solid-solid transition occurs at 
$T=0.16$, while the  population of the liquid-like states passes the  population of   the  icosahedral ones at  $T=0.27$. 
The critical temperature $0.16$ for this network model with harmonic pairwise rates (Eq. \eqref{har}) 
was also obtained in Ref. \cite{cve}   as the point where the transition rates FCC $\rightarrow$ ICO and ICO $\rightarrow$ FCC cross.
This discrepancy in critical temperatures illustrates and quantifies the difference between the continuous $\lj38$ cluster and its network approximation.

\subsection{\label{sec:res}Results}
Originally, the results involving the zero-temperature asymptotics  for the eigenvalues and 
eigenvectors of the $\lj38$ network were presented in Ref. \cite{cam2}. 
Here we highlight the main one of them.
The results associated with the finite temperature continuation are new.
 
\subsubsection{\label{sec:eval}Is there a spectral gap?}
\begin{figure}
\centerline{
\includegraphics[width=0.6\textwidth]{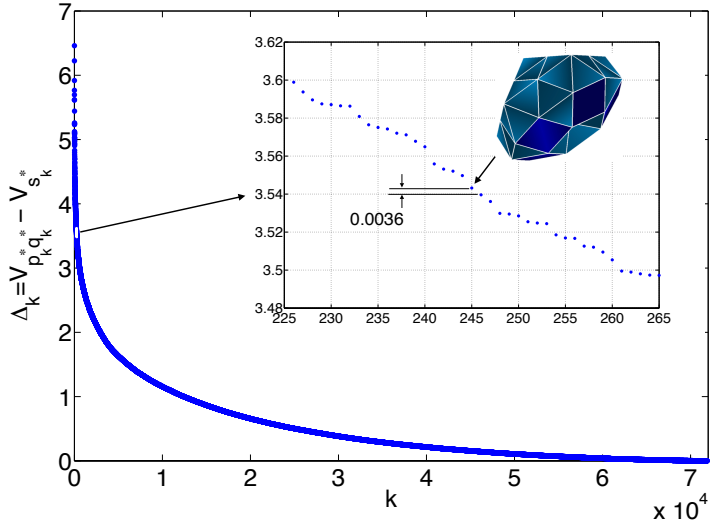}
}
\caption{\label{fig:delta} 
The numbers $\Delta_k$, $k=1,2,\ldots,71886$, 
defining the zero-temperature asymptotics for the eigenvalues of the $\lj38$ network by
$\lambda_k\asymp{-\Delta_k/T}$. The white spot on the graph is blown up in the separate window.
The escape process from the quasi-invariant 
subset of states where the deepest minimum is the second lowest minimum 
ICO corresponds to $k=245$.
Reproduced with permission from NHM 2014. Copyright 2014 American Institute of Mathematical Sciences.
}
\end{figure}

The numbers $\Delta_k$ determining the zero-temperature asymptotics for eigenvalues of the generator 
matrix $L$ of the $\lj38$ network by $\lambda_k\asymp e^{-\Delta_k/T}$ have been originally presented in Ref. \cite{cam2}.
Here we reproduce them for reader's convenience (Fig. \ref{fig:delta}).
It is surprising at the first glance that the numbers $\Delta_k$ 
are arranged pretty densely. 
A few largest $\Delta_k$'s separated by visible gaps 
correspond to high-lying isolated liquid-like local minima 76525 ($k=1$), 91143 $(k=2$), 16497 ($k=3$), etc., 
of no particular interest.
ICO, the second lowest minimum,
corresponds to $k=245$.
The numbers $\Delta_k$ around $k=245$ are zoomed in  in a separate window in Fig. \ref{fig:delta}. The difference between $\Delta_{245}$
and $\Delta_{246}$ is only 0.0036. Therefore, the corresponding eigenvalues 
are separated by an appreciable gap only for extremely low temperatures, at least $T<0.0036$.

\begin{figure}
\centerline{
\includegraphics[width=0.6\textwidth]{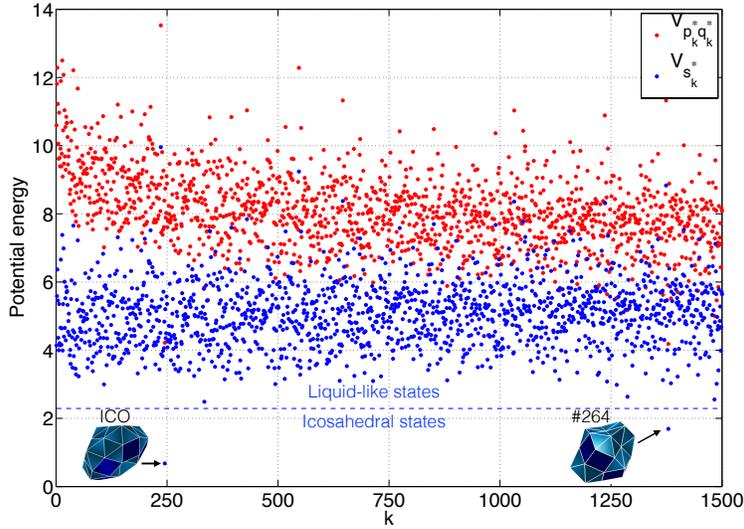}
}
\caption{\label{fig:ms} 
The values  of the potential at the saddles $V_{p^{\ast}q^{\ast}}$ 
and  the local minima $V_{s^{\ast}_k}$ corresponding to the first 1500 sinks.
Among them, only two sinks correspond to icosahedral minima: ICO ($k=245$) and minimum 264 ($k=1379$).
}
\end{figure}

The reason for absence of significant spectral gaps lies in the fact that  
many liquid-like states in the $\lj38$
network are 
$(i)$ separated from the icosahedral and/or the fcc funnels by extremely high saddles
and  $(ii)$ isolated  i.e., correspond to small sets $S_k$ and $C_k$. 
The majority of those sets $S_k$
contain less than 5 states.
The values  of the potential at the saddles $V_{p^{\ast}q^{\ast}}$ and  the local minima $V_{s^{\ast}_k}$ corresponding to the first 1500 sinks 
are shown in Fig. \ref{fig:ms}.
Actually, if, following Ref. \cite{wales_thermo}, we call a state liquid-like if
its energy exceeds $-171.64$ (or 2.29 with respect to FCC), 
the first two states with the smallest positive sink index that are not liquid-like are
ICO  ($k=245$) and state 264 ($k=1379$). 

Roughly speaking, those liquid-like states, separated by extremely high saddles, pollute the interesting part of the 
spectrum, which corresponds to the states that are dynamically accessible from icosahedral or fcc states within 
reasonable waiting times. This fact motivates us to restrict the observation time \cite{wales2014}. This   is  
equivalent to the  introduction of a cap $V_{\max}$ on the potential energy.
Then only the subnetwork $G(S',E')$,
where the sets of states and edges are given by
\begin{align*}
S'&=\{i\in S~\vline~\max_{(p,q)\in w^{\ast}(i,{\rm FCC)}}V_{pq}< V_{\max}\},\\
E'&=\{(i,j)\in E~\vline~i\in S',~j\in S',~V_{ij}< V_{\max}\},
\end{align*}
is considered. 
For $V_{\max}=6$ with respect to FCC (the lowest possible highest barrier separating ICO and FCC is \cite{wales38} 4.219), 
 the resulting network contains 30520 states and 71750 edges.
In this truncated network,
 $\Delta_1$, corresponding to the sink $s^{\ast}_1 = $ ICO, 
is separated by a gap of  approximately 0.19 from $\Delta_2$  
corresponding to  the sink  $s^{\ast}_2 = 223$ (minimum 223) and  the cutting edge with potential energy $V_{p^{\ast}_2q^{\ast}_2}=5.840$.
If we pick $V_{\max}=5.5$, the gap between $\Delta_1$ and $\Delta_2$ will increase to 1.00, while $V_{\max}=5$ will 
make this gap  1.05. 
More details  can be found in Ref. \cite{cam2}.

\subsubsection{Finite temperature eigenvalues, eigenvectors, and eigencurrents}
\begin{figure}
\centerline{
\includegraphics[width=0.6\textwidth]{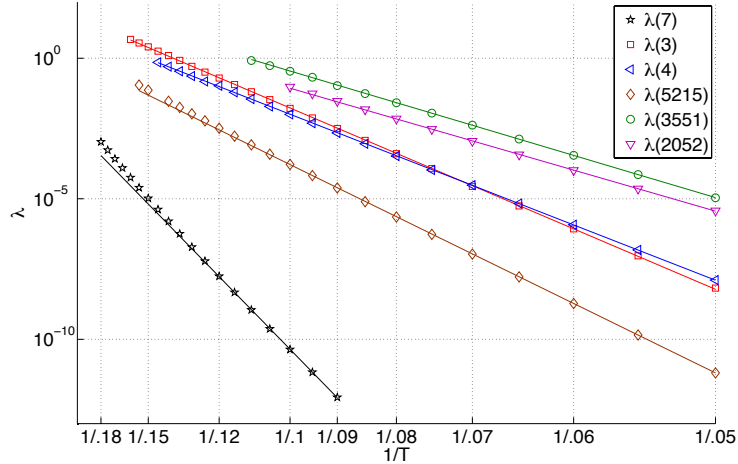}
}
\caption{\label{fig:lam} 
Eigenvalues corresponding to a collection of  large sets $S_k$
are plotted versus $1/T$ where $T$ is  the temperature.
The sinks for these sets are  ICO (minimum 7) and minima  3, 4, 5215, 3551, and 2052.
The lines are obtained by the least squares fit for the low temperature data with linear functions.
}
\end{figure}

\begin{figure}
\centerline{
\includegraphics[width=0.6\textwidth]{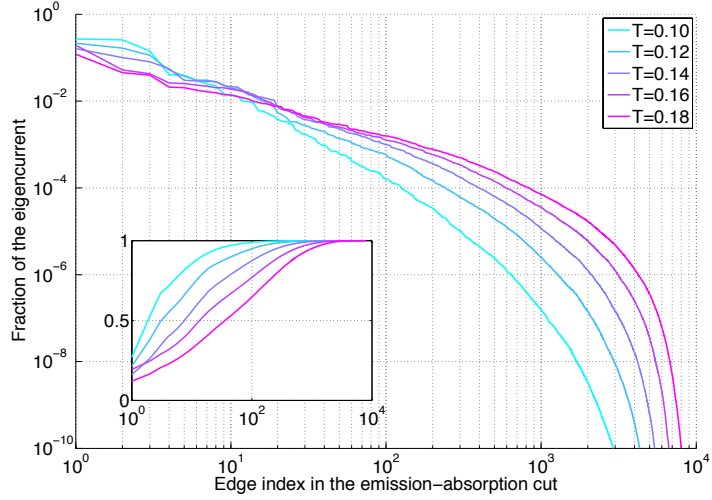}
}
\caption{\label{fig:EAcut} 
The distribution of the eigencurrent in the emission-absorption cut at different temperatures is shown.
The empirical cumulative distribution function is shown in the inset.
}
\end{figure}

The eigenvalues  corresponding to the eigenpairs with the largest sets $S_k$
are plotted in Fig. \ref{fig:lam}. These eigenpairs are identified by \cite{cam2} $k=245$ (sink ICO (7), $|S({\rm ICO})|=56290$),
$k=6910$ (sink 3, $|S(3)|=4252$), $k=7482$ (sink 4, $|S(4)|=1316$), $k=4143$ (sink 5215, $|S(5215)|=990$),
$k=11750$ (sink 3551, $|S(3551)|=680$), and $k=11961$ (sink 2052, $|S(2051)|=1758$).
 The linear leasts squares fits for these eigenvalues  and the numbers $\Delta_k$ are
\begin{equation*}
\begin{array}{lll}
{\rm ICO}: & \Delta = 3.543,&\lambda\approx 1.417\cdot 10^5\cdot e^{-3.570/T}\\
3:& \Delta = 1.408,&\lambda \approx 4.811\cdot 10^4\cdot e^{-1.485/T}\\
4:& \Delta = 1.356,&\lambda \approx 8.222\cdot 10^3\cdot e^{-1.360/T}\\
5215:& \Delta = 1.753,&\lambda \approx 4.253\cdot 10^3\cdot e^{-1.706/T}\\
3551:& \Delta = 1.039,&\lambda \approx 1.122\cdot 10^4\cdot e^{-1.037/T}\\
2052:& \Delta = 1.027,&\lambda \approx 2.029\cdot 10^3\cdot e^{-1.007/T}
\end{array}
\end{equation*}
The exponents of the least squares are close to their asymptotic values. 
Note that the eigenvalues corresponding to sinks 3 and 4 cross at $T=0.07$. 
Also note that the eigenvalue corresponding to sink ICO grows faster than predicted by the Arrhenius law.
This is consistent with the graph of the transition rate from icosahedral 
to fcc states in Ref. \cite{wales_book}. 
This phenomenon is due to the fact that weights of some states other than ICO (e.g. 16 and 8) 
in the left eigenvector $P(T)\phi^{({\rm ICO})}$  grow dramatically as 
the temperature increases. 
As a result, the effective potential at the starting position effectively  increases.
On the other hand, there is an abundance of pathways in
the $\lj38$ network connecting   icosahedral  and fcc states that pass through higher than the lowest possible one but not 
very high saddles \cite{cve}.
This, in turn, leads to a dramatic broadening of the eigencurrent distribution 
in the emission-absorption cut (Fig. \ref{fig:EAcut}). Nearly a power law distribution is observed. Fig. \ref{fig:EAcut}
can be compared to Fig. 9 in Ref. \cite{cve} where the reaction flux distribution in 
the isocommittor cut corresponding to the committor value $q=0.5$ is shown.

\begin{figure}
\centerline{
\includegraphics[width=0.6\textwidth]{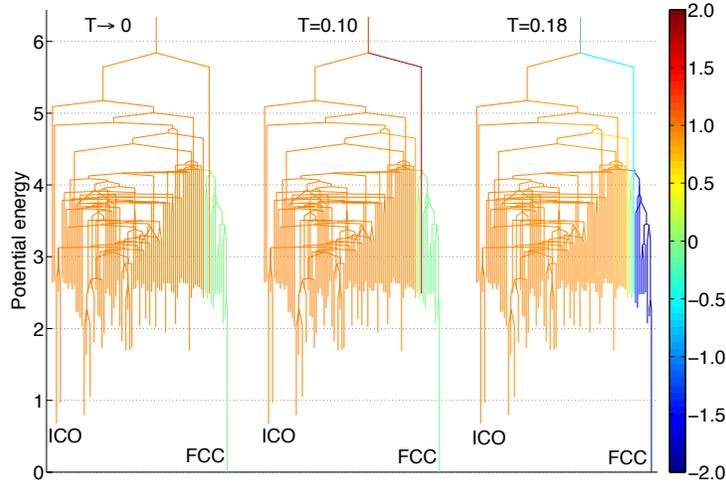}
}
\caption{\label{fig:eve} 
The asymptotic eigenvector (left) and the eigenvectors at $T=0.10$ (middle) and at $T=0.18$ (right) 
corresponding to the sink ICO are 
visualized in the form of disconnectivity graphs. The shown states are selected as follows. 
First, 100 states with the smallest potential energy are selected. Then those that are separated 
by barriers less than 0.2 are lumped together. 
The states are colored according to $\phi^{({\rm ICO})}_j(T)$, $T=0$, 0.10, and 0.18 respectively.
The states are arranged along the $x$-axis in the decreasing order 
according to $\phi^{({\rm ICO})}_j(T=0.18)$.
}
\end{figure}

The asymptotic eigenvector $\phi^{({\rm ICO})}(0)$ corresponding to the sink ICO and its finite temperature  counterparts
$\phi^{({\rm ICO})}(0.10)$ and $\phi^{({\rm ICO})}(0.18)$ at $T=0.10$ and $T=0.18$  respectively are shown in 
Fig. \ref{fig:eve} in the form of disconnectivity graphs. 
There is a close agreement between $\phi^{({\rm ICO})}(0)$ and  $\phi^{({\rm ICO})}(0.10)$,
while the components of $\phi^{({\rm ICO})}(0.18)$  corresponding to the fcc  states  grow  in the negative direction.
This reflects the fact that  icosahedral states are more populated than fcc states in the equilibrium distribution at $T=0.18$.
The majority of the components of  $\phi^{({\rm ICO})}$ remain in the shown  interval $[-2,2]$.
At the same time, a few components corresponding 
to high-lying states  grow large in the absolute value. 
We could not show them without making the colormap non-informative.
The ranges of values of 
 $\phi^{({\rm ICO})}(0.10)$ and   $\phi^{({\rm ICO})}(0.18)$ are $[-9408.0, 577.1]$ and $[-116.8,874.0]$ respectively.
At $T=0.10$, 104 out of 71887 components of $\phi^{({\rm ICO})}$ are greater than 2, and 60 components are less than -2, 
while at $T=0.18$, there are 498 components greater than 2, and  323 components less than -2. 
The emission-absorption cut changes as the temperature grows. In particular, one can spot a state in Fig. \ref{fig:eve}
separated by a high barrier from both  ICO and FCC, that switches sign as the temperature increases from 0.10 to 0.18.

\begin{figure*}
\centerline{
(a)\includegraphics[width=0.75\textwidth]{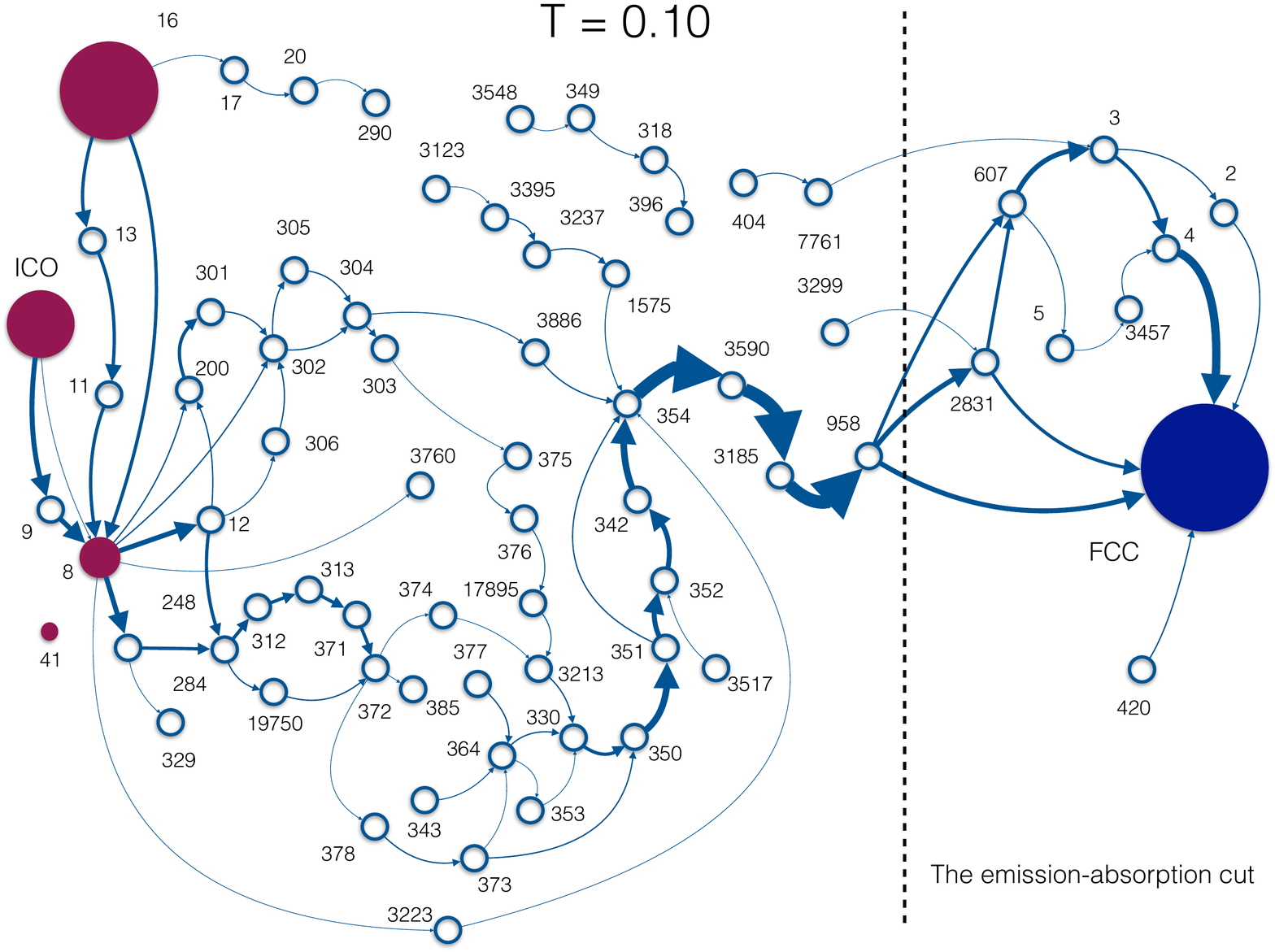}
}\centerline{
(b)\includegraphics[width=0.75\textwidth]{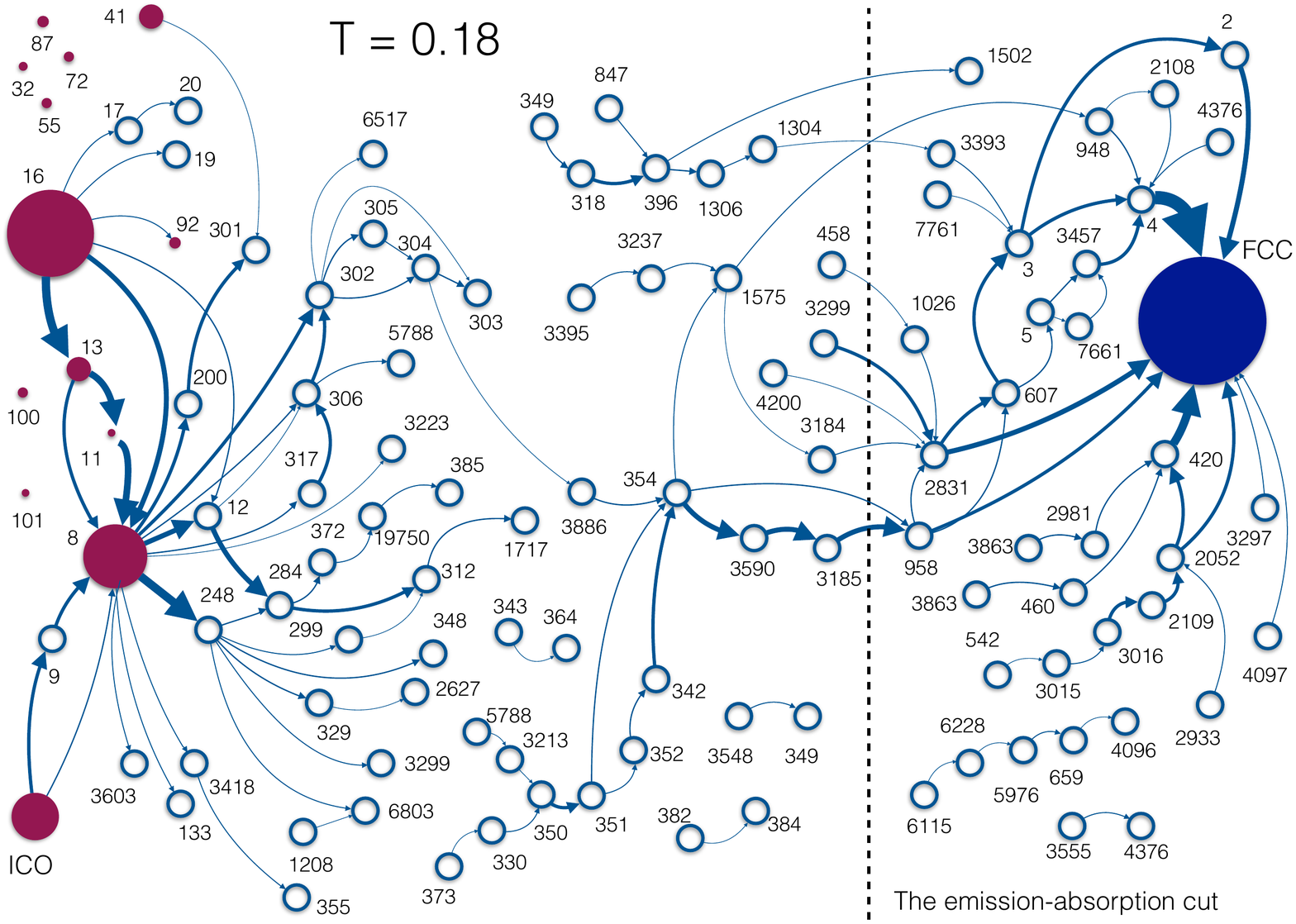}
}
\caption{\label{fig:ec} 
The eigencurrents corresponding to sink ICO at temperatures  $T=0.10$ (a), and $T=0.18$ (b).
The widths of the arrows are proportional to $F^{({\rm ICO})}_{ij}$.
The states which emit or absorb at least 5\% of the eigencurrent are shown by  maroon or bark blue  
solid circles respectively, and the areas of these circles are proportional
to the percentages of emitted/absorbed eigencurrent. 
Other states are represented by empty circles.
The dashed lines indicate the emission-absorption cuts. 
}
\end{figure*}

The relaxation process corresponding to the $k$th eigenpair is quantified by the $k$th eigencurrent defined by Eq. \eqref{eq4}. 
The eigencurrents for  the sink ICO  at temperatures $T=0.10$ and $T=0.18$  are depicted in Fig. \ref{fig:ec} (a) and (b) 
respectively. 
The arrows representing the maximal eigencurrent, 
i.e., $F_{\max}^{({\rm ICO})}(T):=\max_{ij}F^{({\rm ICO})}_{ij}(T)$, have the unit  (the maximal) width.
Only the edges $(i,j)$ with $F_{ij}^{({\rm ICO})}\ge 0.05\cdot F_{\max}^{({\rm ICO})}(T)$ are shown, and their widths 
are proportional to $F_{ij}^{({\rm ICO})}$. 

In Fig. \ref{fig:ec}(a) for $T=0.10$, the sequence of states  
$$354\rightarrow3590\rightarrow3185\rightarrow958$$
stands out by thicker arrows connecting them.
This sequence is a 
part of the asymptotic zero-temperature path (a. k. a. the MinMax path) 
connecting ICO and FCC \cite{cam1}: $\{$ICO, 9, 8, 248, 284, 312, 313, 371, 372, 374, 375, 376, 
17895, 3213, 350, 351, 352, 342, 354, 3590, 3185, 958, 2831, 607, 3, 4, FCC$\}$. 
This path can be traced almost fully in Fig. \ref{fig:ec}(a).
The two highest saddles in the MinMax path are $V_{342,354}=4.219$ and $V_{958,2831}=4.186$.
The height of the saddle between minima 958 and 1, corresponding to another thick arrow across the
emission-absorption cut, is \cite{wales_network} $V_{958,1}=4.272$.

At $T=0.18$, the eigencurrent is considerably more spread  out, 
especially to the right from the emission-absorption cut (see Fig. \ref{fig:ec}(b)).
The number of pathways carrying notable amount of  the eigencurrent from the emission-absorption cut to FCC
dramatically increases, and this cut itself shifts away from FCC. The fraction of the  eigencurrent carried by the sequence
$354\rightarrow3590\rightarrow3185\rightarrow958$ notably drops in comparison with it at $T=0.10$.

The collection of saddles corresponding to the edges of the emission-absorption cut is the natural transition state 
of the relaxation process. 
As shown in Fig. \ref{fig:EAcut}, the distribution of current in this cut is wide, nearly a heavy tailed power law.
This means that the relaxation process dramatically diversifies as the temperature increases.

Figs. \ref{fig:ec} (a) and (b) also show the distributions of the emission and the absorption of the eigencurrent. 
At the range of temperatures $0.09\le T\le 0.18$, minimum 16 carries the maximal weight in the perturbation given by
the left eigenvector $P\phi^{({\rm ICO})}$. The weight  of minimum 8 considerably increases.
At the same time, the weight of ICO, the deepest minimum among the emitting states, drops due to 
its small prefactor (which, in turn, is due to its high order of the point group ($O_{{\rm ICO}} = 10$)).   As the temperature increases, more and more states 
emit appreciable amounts of the eigencurrent.  On the other hand, at this temperature range, FCC remains the dominant absorbing state:
it absorbs over 99\% of the emitted eigencurrent even at $T=0.18$.


\section{\label{sec:tpt}The spectral analysis versus TPT}
In Ref. \cite{cve}, the $\lj38$ network was analyzed using tools of the Transition Path Theory  (TPT) \cite{eve1,eve2,dtpt,luna}.
Throughout this paper we have indicated some points of comparison between the proposed spectral approach and 
TPT. Now we will summarize similarities and differences between these two approaches. 

The basic steps of the analysis  of  a transition process taking place in the 
stochastic network $G(S,E,L)$  by means of TPT are the following.
We assume detailed balance for simplicity.
\begin{enumerate}
\item Pick two disjoint subsets of states: a source set $A\subset S$  and  a sink  set $B\subset S$.
\item Solve the committor equation 
\begin{equation}
\label{com}
\sum_{j\in S\backslash(A\cup B)}L_{ij}q_j=0,\quad q(A)=0,\quad q(B)=1.
\end{equation}
For each state $i$ the value of the committor $q_i$ is the probability that the random walk starting at $i$ will reach first $B$ rather than $A$.
\item Calculate the reactive current
$$F^{R}_{ij}:=\pi_iL_{ij}(q_j-q_i)$$
and analyze it. 
For example, one can consider its distribution in the isocommittor cuts  and visualize the reactive tube with their aid. 
(For a given $q\in[0,1)$ the isocommittor cut \cite{cve} $C(q)$ is the collection of edges $(i,j)$ such that $q_i\le q$ and $q_j>q$.)
\item
Calculate the transition rate $\nu_R$, i.e.,  the average number of transitions from $A$ to $B$ per unit time, e.g. by 
$$\nu_R=\sum_{(i,j)\in C(q)}F^{R}_{ij},$$
where $C(q)$ is an arbitrary isocommittor cut.
Also calculate the rates $k_{A,B}$ and $k_{B,A}$ which are the inverse average times of going back to $B$ after  hitting $A$ and 
going back  to $A$ after hitting $B$ respectively. They are given by
$$k_{A,B}=\frac{\nu_R}{\sum_{i\in S}\pi_i(1-q_i)},\quad  k_{A,B}=\frac{\nu_R}{\sum_{i\in S}\pi_iq_i}. $$
\end{enumerate}

In TPT, the probability distribution in the network is assumed to be equilibrium, i.e., $\pi$.
The transition process from $A$ to $B$ is stationary.
The committor, the reactive current, and the transition rates are time-independent. 
Contrary to this, the eigenpair $(\lambda_k,\phi^k)$ describes the time-dependent relaxation process starting from the 
non-equilibrium distribution $\pi+\alpha P\phi^k$, where the perturbation $\alpha P\phi^k$ decays uniformly throughout the network
with the rate $\lambda_k$. 

Nevertheless, it is instructive to compare the right eigenvector $\phi^k$ to the committor, 
 the eigencurrent 
$$F^k_{ij}=\pi_iL_{ij}(\phi^k_j-\phi^k_i)e^{-\lambda_k/T}$$ 
to the reactive current $F^{R}_{ij}$, and the transition rate $k_{A,B}$ to $\lambda_k$.
We will do it on the example of the $\lj38$ network. 

As the temperature tends to zero, the right eigenvector $\phi^k$ tends \cite{bovier2002,bovier1} to the committor $q^{AB}$
with the source set $A=\{s^{\ast}_0,s^{\ast}_1,\ldots,s^{\ast}_{k-1}\}$ and the  sink set $B=\{s^{\ast}_k\}$, which, in turn, tends to
the indicator function of the set $S_k$.  In Ref. \cite{cve}, the source and sink sets were chosen to be $A\equiv$ ICO and  $B\equiv$ FCC.
Therefore, if we narrow our consideration to the subset of states $S\backslash\bigcup_{j=1}^{244}S_j$ (recall that ICO is the sink $s^{\ast}_{245}$),
we can expect that the right eigenvector $\phi^{({\rm ICO})}$ is close to $1-q^{({\rm ICO,FCC})}$ at low temperatures.
Fig. \ref{fig:eve} visualizing the eigenvector $\phi^{({\rm ICO})}$ at  $T\rightarrow 0$,  $T=0.10$ and $T=0.18$ 
can be compared to Fig. 4 in Ref. \cite{cve} depicting the committor at $T=0.06$, 0.12, and 0.18.

The committor and the right eigenvector play similar roles in the definitions of the corresponding currents.
The net average numbers of transitions per unit time along the edge $(i,j)$
in the Transition Path Process \cite{luna,cve} in TPT and the relaxation process
are proportional to  $\pi_iL_{ij}(q_j-q_i)$ and $\pi_iL_{ij}(\phi^k_j-\phi^k_i)$ respectively.
However, there are important differences between the right eigenvector and the committor. 
While the committor indicates the probability to reach first $B$ rather than $A$ starting from the given state, the
right eigenvector indicates how the given state is over- or underpopulated relative to the equilibrium distribution.
The committor takes values in the interval $[0,1]$, while  the range of  values of the components of the right eigenvector $\phi^{({\rm ICO})}$, 
normalized as described in Section \ref{sec:eig},
is temperature dependent. Some  components acquire values much larger than 1 or much less than -1. 

Both the reactive current $F^R_{ij}$ and the eigencurrent $F^{({\rm ICO})}_{ij}$ describe the net probability flow
for the corresponding processes.
The reactive current is time-independent, 
while the eigencurrent uniformly decays with time at the rate 
$\lambda({\rm ICO})$.
The reactive current is emitted by the source states $A$, absorbed by the sink states $B$ and conserved at all other states \cite{dtpt}.
There is no reactive current between any two source states, and there is no reactive current between any two sink states.
Contrary to this, the eigencurrent is emitted at all states where $\phi^{k}>0$, absorbed at all states where $\phi^k<0$,
and there is a nonzero eigencurrent along any edge $(i,j)$ as long as $\phi^k_i\neq \phi^k_j$. 
The total flux of the reactive current is the same through any cut of the network separating the sets $A$ and $B$.
The total flux of the eigencurrent is maximal through the emission-absorption cut, 
which is the cut separating the states with $\phi^k_j\ge0$ and $\phi^k_j<0$.
Note if $\phi^k_j=0$ at some state $j$, then the emission-absorption cut is 
not the only cut separating the states with $\phi^k>0$ and $\phi^k<0$.
In this case, the flux of the eigencurrent is maximal through any cut separating
the states with $\phi^k>0$ and $\phi^k<0$. 

Despite the differences, the reactive current can serve as a good approximation to 
the quantity $\pi_iL_{ij}(\phi^k_j-\phi^k_i)$ if almost all  eigencurrent is emitted and absorbed by a few states.
For example, this is the case  for the eigencurrent $F^{({\rm ICO})}_{ij}$ in  the $\lj38$ network at low temperatures, 
say $T\le 0.12$.
The distributions of the eigencurrent $F^{({\rm ICO})}_{ij}$ in the 
emission-absorption cut and the reactive current $F^R_{ij}$ in the isocommittor cut corresponding to $q=0.5$ 
are similar (compare Fig. \ref{fig:EAcut} and Fig. 9 in Ref. \cite{cve}).

At low temperatures, both the transition rate $k_{{\rm ICO,FCC}}$ from ICO to FCC 
and the escape rate $\lambda({\rm ICO})$ from the set $S({\rm ICO})$ are well-approximated by the Arrhenius law
with the same parameters. The least squares fits to the Arrhenius law obtained here and in Ref. \cite{cve} are close:
\begin{align*}
k_{{\rm ICO,FCC}}&=9.81\cdot 10^4\cdot e^{-3.525/T},\\
\lambda({\rm ICO})&=1.417\cdot 10^5\cdot e^{-3.570/T}.
\end{align*}
The theoretical asymptotic value of the exponent is $\Delta=3.543$.
As the temperature increases,  the graphs of $k_{{\rm ICO,FCC}}$ and $\lambda({\rm ICO})$ diverge 
(compare Fig. \ref{fig:lam} with Fig. 5 in Ref. \cite{cve}). $\lambda({\rm ICO})$ grows faster than
$k_{{\rm ICO,FCC}}$ due to the fact that the emission of the eigencurrent is distributed.
The fraction emitted by ICO drops  in favor of higher lying states  as the temperature increases.

Finally, we would like to make a few remarks about technical difficulties in applying 
the TPT approach and the proposed spectral approach.
The major computational challenge 
in the TPT lies in solving the committor equation \eqref{com} which is a large linear system. 
However, for networks representing energy landscapes, 
the matrix of the system can be made symmetric. 
In addition, it is positive-definite by construction. Hence,  
the powerful Conjugate Gradient (CG) method with the incomplete Cholesky
preconditioning \cite{demmel} can be used to solve it. 
The main technical difficulty in the proposed spectral approach is solving the large linear system
in the Rayleigh Quotent iteration. The system is symmetric but indefinite,
and its numerical solution is harder than the one of the committor equation. 
Note that the computation of the zero-temperature asymptotics in the spectral approach 
is robust as soon as the genericness assumption
mostly holds. Its violations might lead only to  local effects.


\section{\label{sec:met}Is  LJ$_\mathbf{38}$ metastable?}
A state of an isolated physical system is called metastable, 
if it is not its ground state, but the system can remain in it for a long time.
This concept is clear when the system has just two states, 
but even then some reference time is needed 
in order to determine what period of time is long. 
For example, it can be the observation time. 
It is a nontrivial task to extend the concept of metastability to complex systems.
In a number of works \cite{bovier2002,bovier1,kurchan1,gaveau} the concept of metastability was associated 
with the presence of a spectral gap separating a group of small eigenvalues from the rest of the spectrum.

In this Section, we will discuss whether $\lj38$ is metastable from two points of view: 
of a mathematician and of a chemical physicist. 

A detailed discussion on the metastability of the $\lj38$ network from the point of view of a mathematician can be found in Ref. \cite{cam2}. 
Here we highlight its conclusions.
A formal definition of metastability for the case of Markov chains with detailed balance 
was introduced by Bovier, Eckhoff, Gayrard and Klein \cite{bovier2002,bovier1}.
They have defined the set of metastable points  $\mathcal{M}\subset S$, each of which, in essence, is a  representative of  a 
quasi-invariant subset of states. 
Then the stochastic network is metastable with respect to the set $\mathcal{M}$,
if the  minimal expected time to reach any metastable point from another metastable point is much larger that the maximal
expected time to reach any metastable point from any state not belonging to $\mathcal{M}$.
The absence of the spectral gap in the full $\lj38$ network renders it not metastable  in the sense of Bovier et al 
with respect to a set of metastable points one of which is ICO,
unless the temperature is extremely small,
i.e., $T<0.0036$.  However, capping the potential we can make ICO a metastable point for the considered range of temperatures $0<T\le 0.18$.
Huisinga, Meyn and Schuette proposed another definition of metastability in the context of continuous systems 
where  the time reversibility (the detailed balance)
was not assumed \cite{schuette03,schuette04}. 
Their definition links metastability with ergodicity. 
In their sense, the Freidlin's cycle $C({\rm ICO})$
is metastable for the whole considered range of temperatures $0<T\le0.18$.  

Now we assume the point of view of a chemical physicist.
The Lennard-Jones-38 network is a model for the low temperature dynamics of rare gas clusters of 38 atoms such as argon, krypton and xenon. Rare gas clusters 
have been studied via electron diffraction 
\cite{farges,echt1,harris1,harris2,echt2,kovalenko1,kovalenko2} and extended x-ray absorption fine structure spectroscopy \cite{kakar} experiments.
Mass spectra revealing increased stability of clusters of ``magic" numbers of atoms corresponding to completion (13, 55, 147, 309, ...)
or partial completion (19, 25, ... ) of icosahedral layers were obtained \cite{echt1,harris1,harris2,echt2}. 
The number 38 is not in these lists.
The comparison of the experimentally determined 
mean radial distribution function with the simulated one allowed Kakar et al \cite{kakar} to infer the structure of   clusters for different sizes.
Smaller clusters tend to be icosahedral, while larger ones exhibit signs of the fcc packing. 
The structural transition from icosahedral packing to fcc 
occurs around  cluster sizes from \cite{kakar} $N=200$  to  \cite{farges,harris2} $N=800$ atoms. It has been hypothesized by van der Waal \cite{waal} and confirmed experimentally 
by Kovalenko et al \cite{kovalenko1,kovalenko2} 
that the structural transition from icosahedral packing to fcc occurs via the build up of faulty face-centered cubic layers around
the icosahedral core rather than via cluster rearrangement. Simulated diffraction patterns for faulty fcc clusters were found to be 
in a better agreement with the experimental ones than the patterns simulated from other possible packings. 

These experimental results suggest that the structure of the clusters of 38 atoms generated in experiments is icosahedral.
Otherwise, a peak at 38 atoms in the mass spectra would be present.
Before the structural transition to the global minimum (FCC) has chance to occur, these clusters are likely to acquire more atoms. 
Therefore, the  Lennard-Jones-38 clusters are more likely to be found to be in the icosahedral state that is metastable at low temperatures. 

In order to make our  results comparable with experimental ones we  
divide the collection of states into icosahedral and fcc. 
This partition roughly corresponds to the division into emitting and absorbing states for the
eigencurrent $F^{{\rm ICO}}$. Using parameters appropriate for argon \cite{wales_book} we 
obtain the escape rates listed in Table \ref{tab1}. $T,~*$ and $r,~*$ denote the temperature and the escape rate in the reduced units
used throughout the text.
$T$, K, is the temperature in Kelvins. $r$, s$^{-1}$, is the escape rate in  inverse seconds.
The calculated rates are rather small but comparable with the time scales of argon clusters dynamics \cite{haberland}.
Therefore, if there would be set up an experiment where argon clusters of 38 atoms would be selected and isolated, 
structural rearrangements from icosahedral packing to fcc could be observed. Then the measured rates of such
structural rearrangements could be compared with the calculated rates in Table \ref{tab1}.

\begin{table}
\caption{\label{tab1}Escape rates from icosahedral states to fcc in the $\lj38$ network for parameters appropriate for argon.
The symbol ``*" denotes the reduced units used throughout the text.}
\begin{center}
\begin{tabular}{r|r|l|l|}
$T,~*$ & $T$, K &$r,~*$ &$r$, s$^{-1}$ \\
\hline\\
0.06 & 7.2 & $2.0\cdot 10^{-21}$ & $9.5\cdot 10^{-10}$\\
0.08 & 9.6 & $5.9\cdot 10^{-15}$ & $2.7\cdot 10^{-3}$\\
0.10 & 12.0 & $4.3\cdot 10^{-11}$ & $2.0\cdot 10^{1}$\\
0.12 & 14.4 & $1.8\cdot 10^{-8}$ & $8.1\cdot 10^{3}$\\
0.14 & 16.8 & $1.6\cdot 10^{-6}$ & $7.3\cdot 10^{5}$\\
0.16 & 19.2 & $5.7\cdot 10^{-5}$ & $2.6\cdot 10^{7}$\\
0.18 & 21.6 & $1.1\cdot 10^{-3}$ & $4.9\cdot 10^{8}$\\
\end{tabular}
\end{center}
\end{table}

\section{\label{sec:con}Conclusion}
We have proposed an approach for computing eigenvalues,  eigenvectors, and eigencurrents
of stochastic networks representing energy landscapes.
Step one of this approach is  finding zero-temperature asymptotics for eigenvalues and eigenvectors
using the earlier introduced \cite{cam2}  efficient graph-cutting algorithm.
Step two is the finite temperature continuation using the Rayleigh Quotient iteration supplemented with some precautions.

The proposed methodology has been applied to the $\lj38$ network \cite{wales_network} with 71887 states and 119853 edges.
The results turned out to be stereotype-breaking due to the absence of any significant spectral gaps.
This happens due to the presence of a large number of high-lying liquid-like local minima 
separated by extremely high (perhaps artificial) potential barriers from
the icosahedral and fcc funnels. However, the one-to-one correspondence between local minima and eigenvalues 
obtained by our algorithm renders the task of identifying important eigenvalues trivial.

We have explained the physical meaning of the eigencurrent and proved some of its properties.
We demonstrated how it can be used for quantification and visualization of relaxation processes 
on the example of the $\lj38$ network. Each eigencurrent partitions the set of local minima into emitting and absorbing.
The emission-absorption cut is the true transition state for the relaxation process described by the given  eigencurrent.

We have compared  the results of the analysis of the $\lj38$ network 
by means of the proposed spectral approach and
by means of  TPT. The outcomes of these two analyses are shown to be consistent.

Finally, we have discussed the question of metastability of $\lj38$ from the points of view of a mathematician and a chemical physicist.
A connection with the experimental results for rare gas clusters has been made.


\section*{acknowledgments}
I am grateful to Prof. E. Vanden-Eijnden for  inspiring discussion on  spectra and eigencurrents
and useful suggestions about the present manuscript.
I thank Prof. D. Wales for providing the data  of the $\lj38$ network and for valuable comments.
This work is partially supported by  DARPA YFA Grant N66001-12-1-4220 and NSF Grant 1217118.

\end{document}